\documentclass[12pt]{amsart}
\usepackage{amsmath,amsthm,amsfonts,amssymb}
\usepackage{graphicx}
\usepackage{mathrsfs}
\usepackage{enumitem}
\usepackage{wrapfig}
\usepackage{etoolbox}
\usepackage{xcolor}
\apptocmd{\sloppy}{\hbadness 10000\relax}{}{}
\apptocmd{\sloppy}{\vbadness 10000\relax}{}{}
\usepackage[pdfpagelabels]{hyperref}
\usepackage[letterpaper,margin=1.1in]{geometry}

\numberwithin{equation}{section}
\theoremstyle{plain}

\newtheorem{theorem}{Theorem}[section]

\newtheorem{corollary}[theorem]{Corollary}

\newtheorem{lemma}[theorem]{Lemma}

\theoremstyle{definition}
\newtheorem{remark}[theorem]{Remark}
\newtheorem{definition}[theorem]{Definition}

\def\XXint#1#2#3{{\setbox0=\hbox{$#1{#2#3}{\int}$ }
\vcenter{\hbox{$#2#3$ }}\kern-.6\wd0}}

\newcommand{\dist}{\mathop\mathrm{dist}\nolimits}

\newcommand{\loc}{\mathrm{loc}}

\usepackage{alphalph}

\title{Non-existence of cusps for a Free-boundary Problem for Water Waves}
\author{Sean McCurdy}

\keywords{Stokes wave, partial regularity}

\address{Department of Mathematics\\ National Taiwan Normal University\\ Taipei, Taiwan}
\email{smccurdy@ntnu.edu.tw}

\begin{document}

\maketitle

\begin{abstract} In \cite{WeissVarvaruca11}, Varvaruca and Weiss eliminate the existence of cusps for a free-boundary problem for two-dimensional water waves under assumptions that hold for solutions such that $\{u>0\}$ is a ``strip-like" domain in the sense of \cite{Varvaruca08}.  In this paper it is proven that cusps do not exists in the natural setting for these free-boundary problems. In particular, non strip-like domains are also allowed.  This builds upon recent work on non-existence of cusps in \cite{MccurdyNaples22}.       
\end{abstract}

\section{Introduction}
In this note, we eliminate the existence of cusps in the free-boundary for weak solutions of the following free-boundary problem
\begin{align} \label{e:free-boundary problem}\tag{$P_{\gamma}$}
\Delta u= 0 & \qquad \text{in  } \{u>0\} \cap \Omega \subset \mathbb{R}^n\\
\nonumber |\nabla u(x_1, ..., x_n)| = |x_n|^{\gamma}  & \qquad \text{on  } \partial \{u>0\} \cap \Omega
\end{align}
for $n=2$ any $0< \gamma$.  For $n=2$, $\gamma=\frac{1}{2}$, these free-boundary problems have a history dating back to 1847 and the work of Stokes \cite{Stokes1880} on $2$-dimensional inviscid, incompressible fluids acted upon by gravity with a free surface.  Stokes studied the profiles of standing waves for such a fluid under for $\Omega = [-1, 1]\times[-D, 1]$, $1<D\le \infty$, and imposing the physical boundary condition that $u$ is constant on $[-1, 1]\times\{-D\}$. Stokes conjectured that there was a one-parameter family of solutions to (\ref{e:free-boundary problem}) parametrized by wave height.  This family of \textit{stream functions} was conjectured to have a wave of maximal height for which the wave profile $\partial \{u>0\}$ touches $\{x_2=0\}$ with angle $\frac{2\pi}{3}$.  In honor of Stokes, this extremal wave is called the Stokes wave.  Points in $\partial \{u>0\} \cap \{x_2=0 \}$ are called \textit{stagnation points}.  Under strong assumptions assumptions of symmetry, monotonicity, and graphicality of $\partial \{u>0\}$, Toland \cite{Toland78} and McLeod \cite{McLeod97} proved the existence of extreme periodic waves, for $R = + \infty$ and $0<R<\infty$ (i.e., for waves of infinite and finite depth). In 1982, Amick, Fraenkel and Toland \cite{AmickFraenkelToland82} and Plotnikov \cite{Plotnikov82} both independently proved the Stokes' conjecture on the aperture of $\{u>0\}$ where it touches $\{x_2=0\}$ for this extremal wave under similar assumptions.  See \cite{AramaLeoni12, WeissVarvaruca11} for more historical details and the derivation of (\ref{e:free-boundary problem}) from physical principles.

Solutions to the free boundary problem (\ref{e:free-boundary problem}) are critical points of the corresponding Alt-Caffarelli functional 
\begin{align}\label{e:functional}
J_{\gamma}(v) := \int_{\Omega}|\nabla v|+x_n^{2\gamma}\chi_{\{u>0 \}}dx.    
\end{align}
Because they are merely critical points and not minimizers, research on the Stokes Wave has typically proceeded by analyzing  \textit{weak solutions} \cite{ShargorodskyToland08, Varvaruca08}.  These weak solutions assume a modicum of regularity to $\partial \{u>0\}$ away from $\{x_n = 0\}$, see Definition \ref{weak solution def}.  However, in recent work \cite{WeissVarvaruca11}, Varvaruca and Weiss introduced the notion of a \textit{variational solution} (see \cite{WeissVarvaruca11} Definition 3.1). This variational notion of a solution allows Varvaruca and Weiss to employ geometric techniques to the study of the Stokes wave beyond the usual assumptions of symmetry, monotonicity, and graphicality of $\partial \{u>0\}$.  These geometric techniques are based upon a monotonicity formula analogous to the quantity in \cite{Weiss99} for local minimizers of the Alt-Caffarelli functional $J_{0}$. In particular, this allowed them to obtain the following results.

\begin{theorem}\label{sigma decomposition}\emph{(\cite{WeissVarvaruca11})}
Let $n \ge 2$, $\gamma = 1/2$, $\Omega \subset \mathbb{R}^n$ be open, and let $u$ be a weak solution of (\ref{e:free-boundary problem}) such that locally $u$ satisfies
\begin{align}\label{gradient bound}
    |\nabla u| \le C_0 x_n^{\frac{1}{2}}.
\end{align}
Then, if we denote $\Sigma := \partial \{u>0 \} \cap \{x_n=0\} \cap \Omega$ we may decompose $\Sigma$ into the disjoint union $\Sigma = \Sigma_{\text{cusps}} \cup \Sigma_{\text{rect}} \cup \Sigma_{\text{iso}},$ where if we denote the density $\theta^n_{\{u>0\}}(x) = \lim_{r\rightarrow 0^+}\frac{\mathcal{H}^n(B_r(x) \cap \{u>0\})}{\omega_n r^n}$, then the density exists and we define
\begin{align*}
\Sigma_{\text{cusps}} & := \left\{x \in \Sigma: \theta^n_{\{u>0\}}(x) =0\right\}\\
\Sigma_{\text{rect}} & := \left\{x \in \Sigma : \theta^n_{\{u>0\}}(x) \in (0, 1/2)\right\}\\  
\Sigma_{\text{iso}} & := \left\{x \in \Sigma : \theta^n_{\{u>0\}}(x)=1/2\right\}.
\end{align*}
Furthermore, these sets satisfy the following properties.
\begin{enumerate}
    \item $\dim_{\mathcal{H}}(\Sigma_{rect}) \le n-2$, and if $n=2$, then $\Sigma_{rect}$ is locally isolated \cite{WeissVarvaruca11}.
    \item $\Sigma_{\text{rect}}$ is countably $(n-2)$-rectifiable, $\overline{\dim}_{\mathcal{M}}(\Sigma_{rect}) \le n-2$, and if $n=2$, then $\Sigma_{rect}$ is locally finite \cite{WeissVarvaruca11} with bound $\mathcal{H}^0(\Sigma_{\text{rect}} \cap B_r(x_1, 0)) \le C(n, C_0)$ if (\ref{gradient bound}) holds in $B_{2r}(x_1, 0)$ \cite{Mccurdy20}.
    \item $\Sigma_{\text{iso}}$ is closed, and if $n=2$ then $\Sigma_{\text{iso}}$ is locally finite \cite{WeissVarvaruca11}.
\end{enumerate}
\end{theorem}

We note that the assumption (\ref{gradient bound}) is absolutely essential to the project of \cite{WeissVarvaruca11}.  Assumption (\ref{gradient bound}) provides the compactness necessary to geometric blow-up analysis.  And, it is also necessary in order to connect weak solutions and variational solutions, since any weak solution of (\ref{e:free-boundary problem}) which also satisfies (\ref{gradient bound}) is a variational solution. See \cite{WeissVarvaruca11} Lemma 3.4 for details.

The main result of this paper is to eliminate the set $\Sigma_{cusps}$.

\begin{theorem}\emph{(Main Result)}\label{VarvarucaWeiss improvement}
Let $n=2$ and $0< \gamma$.  Let $u$ be a weak solution to (\ref{e:free-boundary problem})  such that for all $(x, 0) \in \Omega$ there exists a neighborhood $K \subset \subset \Omega$ of $(x, 0)$ and a constant $C< \infty$ (possibly depending upon $K$) such that
\begin{align}\label{checkable}
    |\nabla u(x_1, x_2)| \le C|x_2|^{\gamma},
\end{align}
for all $(x_1, x_2) \in K$. Then $\Sigma_{\text{cusps}} = \emptyset$.
\end{theorem}

Theorem \ref{VarvarucaWeiss improvement} is inspired by and improves upon the following result from \cite{WeissVarvaruca11}.

\begin{theorem}\label{cusps needs improvement}\emph{[\cite{WeissVarvaruca11}, Lemma 4.4]}
Let $n=2$, $\gamma = \frac{1}{2}$, and let $u$ be a weak solution to \ref{e:free-boundary problem} satisfying 
\begin{align}\label{e: uncheckable}
    |\nabla u(x_1, x_2)| \le  |x_2|^{\frac{1}{2}}.
\end{align}
Then $\Sigma_{\text{cusps}} = \emptyset$.
\end{theorem}

\begin{remark}
The assumption (\ref{e: uncheckable}) is much stronger than (\ref{checkable}).  In particular, it is not known whether or not weak solutions, in general, satisfy (\ref{e: uncheckable}).  For example, local minimizers of the corresponding Alt-Caffarelli functional $J_{\gamma}$ are weak solutions to (\ref{e:free-boundary problem}) and satisfy (\ref{checkable}) for a dimensional constant $1< C= C(n, \gamma)$ \cite{Mccurdy20}. One might expect weak solutions, which are merely critical points, to be less well-behaved than local minimizers.  In \cite{Varvaruca08} it is proven that if $u$ is a solution to (\ref{e:free-boundary problem}) for any $0< \gamma$ and  $\{u>0\}$ is a ``strip-like" domain (see Section 2.1 in \cite{Varvaruca08}), then $u$ satisfies \eqref{e: uncheckable} (see the proof of \cite{Varvaruca08} Theorem 3.6, in particular it follows from the properties of the function $Q$ in (4.19)). However, ``strip-like" domains do not allow for air bubbles, and therefore Theorem \ref{cusps needs improvement} only represents a partial solution to eliminating $\Sigma_{cusps}$.

The central improvement of this paper is to eliminate the existence of $\Sigma_{cusps}$ under the more natural assumption (\ref{gradient bound}) and only using the local properties a cusp must satisfy.  In particular, this allows for $\{u>0\}$ which are not ``strip-like" in the sense of \cite{Varvaruca08} and for solutions to (\ref{e:free-boundary problem}) which do not satisfy the boundary conditions of wave equations.  

The method of proof for Theorem \ref{VarvarucaWeiss improvement} was initially inspired by the proof of Theorem \ref{cusps needs improvement}. However, the improvement comes from a closer analysis of local cusp geometry using ideas introduced in \cite{MccurdyNaples22}. 
\end{remark}

In fact, Theorem \ref{VarvarucaWeiss improvement} is a qualitative result which comes from a ``quantitative" result on the geometry of the free boundary $\partial \{u>0\}$. To state the result, we need to first define a family of rescalings and a height function, which will be central to helping us control the geometry of $\partial \{u>0\}$ near $(x, 0) \in \partial \{u>0\} \cap \Omega$.

\begin{definition}\label{rescalings}(Rescalings)
Let $0<\gamma$, and let $u$ a weak solution to (\ref{e:free-boundary problem}) in the domain $\Omega$ .  For any set $U \subset \Omega$, $(x, 0)$, and $0<r$ we define the rescalings
\begin{align}
    U_{(x, 0), r} := \frac{U - (x, 0)}{r}.
\end{align} 
If $(x, 0) \in \partial \{u>0\} \cap \Omega$ then
\begin{align}
    u_{(z, 0), r}(x_1, x_2) := \frac{u(r(x_1, x_2)+(x, 0))}{r^{\gamma}}
\end{align}
is a solution to (\ref{e:free-boundary problem}) in $\Omega_{(x, 0), r}$ and $\{u_{(x, 0),r}>0 \} = \{u>0\}_{(x, 0), r}$.

If $(x, 0) \in \partial \{u>0\} \cap \Omega$ then for any $0< r $ there must be a component $\mathcal{O}_{(x, 0), r}$ of $\{u_{(x, 0), r}>0 \} \cap [-1, 1] \times [-1, 1]$ such that $(0, 0) \in \partial \mathcal{O}_{(x, 0), r}$. Furthermore, it is clear that we can choose the components $\mathcal{O}_{(x, 0), r}$ in a consistent manner such that for all $0< r_1< r_2$, $\mathcal{O}_{(x, 0), r_1} = (\mathcal{O}_{(x, 0), r_2})_{(0,0), \frac{r_1}{r_2}}$. 
\end{definition}

\begin{definition}(Height Function)
Let $0< \gamma$, and let $u$ is a weak solution of (\ref{e:free-boundary problem}) such that $(x, 0) \in \partial \{u>0\} \cap \Omega$. Then, for any $0< r \le 2^{-1}\dist((x, 0), \partial \Omega)$ and any component $\mathcal{O}_{(x, 0), r}$ of $\{u_{(x, 0), r}>0 \} \cap [-1, 1] \times [-1, 1]$ such that $(0, 0) \in \overline{\mathcal{O}_{(x, 0), r}}$, we define the following function. For $0< \rho \le 1$, we define
\begin{align*}
    \text{Height}(\rho, \mathcal{O}_{(x,0), r}):= & \min \{ 1, \sup \{|x_2| : (x_1, x_2) \in \mathcal{O}_{(0,0), r}, \quad |x_1| = \rho\} \}.
\end{align*}
Note that $\text{Height}(\rho, \mathcal{O}_{(x,0), r}) = \text{Height}(1, \mathcal{O}_{(x,0), \rho r})= \text{Height}(\rho r, \mathcal{O}_{(x,0), 1})$.
\end{definition}

\begin{theorem}\label{macroscopic information}\emph{(Quantitative Result)}
Let $n=2$ and $0< \gamma$.  Let $u$ be a weak solution to (\ref{e:free-boundary problem}) with associated domain $\Omega$.  Suppose that $(x, 0) \in \Omega \cap \partial \{u>0\}$ and $u$ satisfies (\ref{checkable}) in $B_r((x, 0)) \subset \Omega$ with constant $C< \infty$. Then $\emph{Height}(1, \mathcal{O}_{(x, 0), \rho}) \ge \frac{1}{6C}$ for all $0< \rho \le r$.
\end{theorem}


\begin{remark}
The proofs of Theorem \ref{VarvarucaWeiss improvement}, Theorem \ref{cusps needs improvement}, and Theorem \ref{macroscopic information} are essentially restricted to $n = 2$. For recent results eliminating cusps in higher dimensions and arbitrary co-dimension, see \cite{MccurdyNaples22} which obtained an analogous macroscopic geometric description of $\partial \{u>0 \}$ for local minimizers of an analogous Alt-Caffarelli functional $J_{\gamma}$ in $n \ge 2$.
\end{remark}

It is unknown whether or not $\Sigma_{cusps}$ may be eliminated in $n \ge 3$. It is unknown whether or not $\Sigma_{iso}$ may be eliminated in $n\ge 2$.

\subsection{Acknowledgements}
This author thanks Giovanni Leoni for posing the question addressed in this paper.  More broadly, the author is deeply indebted to Giovanni Leoni and Irene Fonseca for their invaluable generosity, patience, and guidance, and acknowledges the Center for Nonlinear Analysis at Carnegie Mellon University for its invaluable support.

\section{Preliminaries and Reduction of Theorem \ref{VarvarucaWeiss improvement} to Theorem \ref{macroscopic information}}

We begin with by defining the appropriate notion of a solution to (\ref{e:free-boundary problem}).

\begin{definition}\label{weak solution def}\text{(Weak Solutions)}
Let $\Omega \subset \mathbb{R}^n$ and $0< \gamma$.  A function $u \in W^{1,2}_{loc}(\Omega)$ is a \textit{weak solution} of (\ref{e:free-boundary problem}) if $u$ satisfies,
\begin{enumerate}
    \item $u \in C^{0}(\Omega)$, $u \ge 0$ in $\Omega$.
    \item $u$ is harmonic in $\{u>0\} \cap \Omega$.
    \item For every $\tau >0$, the topological free boundary, $\partial \{u>0\} \cap \Omega \cap \{|x_n| >\tau \}$, can be decomposed into an $(n-1)$-dimensional $C^{2, \alpha}$-surface, denoted $\partial_{red}\{u>0\}$, which is relatively open in $\partial \{u>0\},$ and a singular set with $\mathcal{H}^{n-1}$-measure zero.
    \item For any open neighborhood $V$ containing a point $x_0 \in \Omega \cap \{|x_n| > \tau \} \cap \partial_{red}\{u>0\}$, the function $u \in C^1(V \cap \overline{\{u>0\}})$ and satisfies $|\nabla u|^2 = x_n^{2\gamma}$ on $V \cap \partial_{red}\{u>0\}$.
\end{enumerate}
\end{definition}

\begin{remark}\label{r:physical assumptions on weak solutions} 
We note that for physical reasons, the definition of a weak solution also includes the assumption $u \equiv 0$ in $\Omega \cap \{x_n \le 0\}$. However, (\ref{e:free-boundary problem}) is only a physical problem for $n=2$ and $\gamma = 1/2$. In this note, we work without this assumption to allow a wider class of solutions.
\end{remark}.




In order to reduce Theorem \ref{VarvarucaWeiss improvement}
to the proof of Theorem \ref{macroscopic information}, we need to following compactness result.  

\begin{lemma}\label{l:VW 4.1}
Let $n=2$ and $0<\gamma$. Let $u$ be a weak solution which satisfies (\ref{checkable}), $(x, 0) \in \Sigma$, and $0< r_i \rightarrow 0$.  Then, there is a $(1+\gamma)$-homogeneous function $u_{\infty} \in W^{1,2}_{loc}(\mathbb{R}^2)$ and a subsequence such that $u_{(x, 0), r_i} \rightarrow u_{\infty}$ in $W^{1, 2}_{loc}(\mathbb{R}^2)$.
\end{lemma}

\begin{proof}
The case $\gamma = 1/2$ is proven in \cite{WeissVarvaruca11} Lemma 4.1 under the physical assumptions on weak solutions in Remark \ref{r:physical assumptions on weak solutions}.  This argument holds \textit{verbatim} for $0< \gamma< \infty$. To wit, the assumption (\ref{checkable}) implies that for $0<R<\infty$ and all sufficiently large $i$ (depending upon $R$) the functions $u_{(x, 0), r_i}$ are uniformly bounded in $W^{1, 2}(B_R(0))$.  By Rellich-Kondrachov compactness and lower semicontinuity, it remains to show that if $u_{(x, 0), r_i} \rightharpoonup u_{\infty}$ in $W^{1, 2}(B_R(0))$, then $$\int_{\mathbb{R}^2}|\nabla u_{\infty}|^2 \eta d\mathcal{H}^2 \ge \limsup_{i \rightarrow \infty}\int_{\mathbb{R}^2}|\nabla u_{(x, 0), r_i}|^2 \eta d\mathcal{H}^2$$ for all $\eta \in C^1_0(B_R(0))$.  By assumption (\ref{checkable}) and Arzela-Ascoli, $u_{\infty}$ is continuous and $u_{(x, 0), r_i} \rightarrow u_{\infty}$ uniformly in $B_R(0)$. Since $u_{(x, 0), r_i}$ are harmonic in $\{u_{(x, 0), r_i}>0\}$, $u_{\infty}$ is also harmonic in $\{u_{\infty}>0\}$. Therefore, by integration by parts we calculate
\begin{align*}
    \int_{\mathbb{R}^2}|\nabla u_{(x, 0), r_i}|^2 \eta d\mathcal{H}^2 & = - \int_{\mathbb{R}^2} u_{(x, 0), r_i} \nabla u_{(x, 0), r_i} \cdot \nabla \eta d\mathcal{H}^2\\
    & \rightarrow - \int_{\mathbb{R}^2} u_{\infty} \nabla u_{\infty} \cdot \nabla \eta d\mathcal{H}^2 = \int_{\mathbb{R}^2}|\nabla u_{\infty}|^2 \eta d\mathcal{H}^2.
\end{align*}

It remains to show that $u_{\infty}$ is $(1+ \gamma)$-homogeneous.  This is proven in \cite{Mccurdy20} Theorem 5.11 for $n=2$, $k=1,$ and $\Gamma = \{x_2=0\}$ for local minimizers of $J_{\gamma}$.  However, since weak solutions which satisfy (\ref{checkable}) also satisfy the monotonicity formula in \cite{Mccurdy20} Theorem 4.3, Theorem 5.11 holds for them as well.
\end{proof}

\subsection{The reduction of Theorem \ref{VarvarucaWeiss improvement} to Theorem \ref{macroscopic information}}

If we assume Theorem \ref{macroscopic information}, then Theorem \ref{VarvarucaWeiss improvement} will follow if it can be shown that if $(x, 0) \in \Sigma_{cusps},$ then for any $0<C_2<1$, there exists a radius $0<r$ such that $\text{Height}(1, \mathcal{O}_{(x,0), r}) \le C_2$. This follows from an argument analogous to \cite{WeissVarvaruca11} Lemma 4.4. 

To wit, we may consider $\Delta u$ as a non-negative Radon measure supported on $\partial \{u>0\}$ which satisfies the following inequality.
\begin{align*}
    \Delta u(\phi) & := -\int_{\mathbb{R}^2} \nabla u \cdot \nabla \phi d\mathcal{H}^2 \qquad \text{for   } \phi \in C^{\infty}_c(\mathbb{R}^2)\\
    \Delta u(U) & := \sup \{\Delta u(\phi): \phi \in C^{\infty}_c(U), |\phi|_{\infty} \le 1\}\\
    & \ge \int_{\partial_{red}\{u>0\} \cap U} |x_2|^{\gamma}d\sigma(x_1, x_2).
\end{align*}
Now, let $u_{(x, 0), r}$ be the rescaling of $u$, and let $u_{(x, 0), r}|_{\mathcal{O}_{(x,0), r}}$ be the piece-wise function
\begin{align*}
    u_{(x, 0), r}|_{\mathcal{O}_{(x,0), r}} = \begin{cases}
    u_{(x, 0), r} & \text{ in } \mathcal{O}_{(x,0), r}\\
    0 & \text{ elsewhere}.
    \end{cases}
\end{align*}
For $0<r$ small enough, $u_{(x, 0), r}|_{\mathcal{O}_{(x,0), r}}$ is a weak solution in $B_2(0, 0)$ and hence by Lemma \ref{l:VW 4.1}, there exists a function $u_0 \in W_{\loc}^{1, 2}(\mathbb{R}^2)$ such that $u_{(x, 0), r}|_{\mathcal{O}_{(x,0), r}} \rightarrow u_0$  in $W_{loc}^{1, 2}(\mathbb{R}^2)$
as $r \rightarrow 0$.  By the assumption that $\Theta^2((x, 0), \mathcal{O}) = 0$ and the fact that $u_0$ is homogeneous, we have that $u_0 \equiv 0$ and $\Delta u_0 \equiv 0$. 

Thus, $\Delta (u_{(x, 0), r}|_{\mathcal{O}_{(x,0), r}}) \rightharpoonup \Delta u_0$ as $r \rightarrow 0$ and 
\begin{align}\label{flattening}\nonumber
    \Delta u_{(x, 0), r}|_{\mathcal{O}_{(x,0), r}}(B_2^2(0,0)) & \ge  \int_{\partial_{red} \mathcal{O}_{(x,0), r} \cap B_2^2(0,0)} |y|^{\gamma}d\sigma \\
    & \gtrsim \left|\sup_{0< \rho< 1}\{|x_2|: (x_1, x_2) \in \partial \mathcal{O}_{(x, 0), r} \cap B_2^2(0, 0), |x_1|=\rho \}\right|^{2\gamma} \rightarrow 0.
\end{align}
Thus, Theorem \ref{VarvarucaWeiss improvement} follows from Theorem \ref{macroscopic information}.

\section{Proof of Theorem \ref{macroscopic information}}

\subsection{Main Geometric Observations}

The following geometric observation was inspired by the insight that if $(x, 0) \in \Sigma_{cusps}$, then $\partial \{u>0 \}$ must approach $\{x_2 = 0\}$ tangentially in the sense of \eqref{flattening}. That is, $\theta^2((x, 0), \mathcal{O}) = 0$ and (\ref{flattening}) implies that $\text{Height}(1, \mathcal{O}_{(x, 0), r}) \rightarrow 0$ as $r\rightarrow 0$.  Therefore, we may not expect a component of $\{u>0\}$ which touches $(x, 0)$ to be contained in a set of the form $\{ (x_1, x_2) : |x_2| \ge m|x_1|\}$ in any neighborhood of $(x, 0)$ for any $0<m$.  The content of the lemma below is that if we weaken this to consider sets of the form $\{(x_1, x_2):|x_2|\ge m|x_1|-b\}$, then for appropriate choices of $0< m, b$ we can find a neighborhood in which a large piece of $\partial\{u>0\}$ is contained in such a set.

\begin{lemma}\label{n=2 geometry}(Main Geometric Observation)
Let $n=2$, $0< \gamma$, and let $u$ be a weak solution to (\ref{e:free-boundary problem}) that satisfies (\ref{checkable}). Assume that $(x, 0) \in \Sigma_{cusps}$ with $\delta = \dist((x, 0), \partial \Omega)$.  Let $\mathcal{O}$ be a component of $\{u>0\} \cap B_{\delta}((x, 0))$ such that $(x, 0) \in \partial \mathcal{O}$.  For any $0< C_2 \le \frac{1}{2}$ if there exists a radius $0< r_0 \le \delta$ such that 
\begin{align*}\nonumber
    \emph{Height}(1, \mathcal{O}_{(x,0), r_0}) & \le C_2,
\end{align*}
then there exists a $0<\rho \le r_0$ such that $\emph{Height}(1, \mathcal{O}_{(x,0), \rho}) \le C_2$ and one of the following conditions hold
\begin{align}\label{L_r}
    \emph{Height}(r, \mathcal{O}_{(x,0), \rho}) & \ge  3 \emph{Height}(1, \mathcal{O}_{(x,0), \rho})r - 2\emph{Height}(1, \mathcal{O}_{(x,0), \rho})
\end{align}
for all $r \in [2/3, 1]$, or
\begin{align}\label{R_r}
    \emph{Height}(r, \mathcal{O}_{(x,0), \rho}) & \ge  -3 \emph{Height}(1, \mathcal{O}_{(x,0), \rho})r - 2\emph{Height}(1, \mathcal{O}_{(x,0), \rho})
\end{align}
for all $r \in [-1, -2/3]$.
\end{lemma}

\begin{remark}
The lines described by equality in (\ref{L_r}) and (\ref{R_r}) are the lines which intersects the points $(\pm 1, \text{Height}(1, \mathcal{O}_{(x,0), \rho}))$ and $(\pm 2/3, 0)$, respectively.
\end{remark}

\begin{proof}
Let $u$, $0<r_0$, $0< C_2<1/2$ be given.  For ease of notation, we note that by reflection across the $x$- and $y$-axes if necessary, we may assume that $u$ attains $\text{Height}(1, \mathcal{O}_{(x,0), r_0})$ in $\{(1, x_2): x_2 \in \mathbb{R}_+\}$ and not $\{(1, x_2): x_2 \in \mathbb{R}_-\}$ or $\{(-1, x_2) : x_2 \in \mathbb{R}\}$.

Let $L_r$ be the line given by the graph of the function 
\begin{align*}
x_2 & = 3\text{Height}(r, \mathcal{O}_{(x,0), r_0}) x_1 - 2 \text{Height}(r, \mathcal{O}_{(x,0), r_0}).
\end{align*}
We claim that we can find an $0 < \frac{1}{2} <r \le 1 $ such that
\begin{align*}
    \text{Height}(x_1, \mathcal{O}_{(x,0), r_0}) & \ge 3\text{Height}(r, \mathcal{O}_{(x,0), r_0})x_1 - 2 \text{Height}(r, \mathcal{O}_{(x,0), r_0}).
\end{align*}
for all $\frac{2}{3}r \le x_1 \le r$. If we can find such a radius, then $\rho = r \cdot r_0$ proves the lemma.

To prove the claim, we argue by contradiction.  Let $r_1=1$.  If $r_i$ does not satisfy the claim, then there must exist a radius $\frac{2}{3}r_i < r < r_i$ such that 
\begin{align}\label{peaks below}
    \text{Height}(r, \mathcal{O}_{(x,0), r_0}) & < 3\text{Height}(r_i, \mathcal{O}_{(x,0), r_0}) r - 2 \text{Height}(r_i, \mathcal{O}_{(x,0), r_0}).
\end{align}
Let $r_{i+1} \in [2r_i/3, r_i]$ be defined by 
\begin{align*}
    r_{i+1} := \inf\{r \in (2r_i/3, r_i): \text{(\ref{peaks below}) holds}\}.
\end{align*}
Observe that by construction
\begin{align} \label{Height monotonicity}
\text{Height}(r_{i+1}, \mathcal{O}_{(x,0), r_0}) & < \text{Height}(r_i, \mathcal{O}_{(x,0), r_0}).
\end{align}

If the inductively defined sequence $\{r_i \}_i$ does not terminate in finitely many steps with a radius which satisfies the claim, then $\{r_i \}_i$ forms a monotonically decreasing sequence in $[1/2, 1]$, and there is a limit point $r_\infty \in [\frac{1}{2}, 1]$ such that $r_i \rightarrow r_{\infty}$.  By (\ref{Height monotonicity}), there are two possibilities: either $\text{Height}(r_\infty, \mathcal{O}_{(x,0), r_0}) > 0$ or $\text{Height}(r_\infty, \mathcal{O}_{(x,0), r_0}) = 0$. The latter case contradicts the assumption that $\mathcal{O}_{(x,0), r_0}$ is a connected component which touches $(0,0)$. Therefore, we may assume that $\text{Height}(r_{\infty}, \mathcal{O}_{(x,0), r_0})>0$.  We claim that $r_{\infty} = r$. 

By the convergence of $\{r_i\}$ and the fact that $\{\text{Height}(r_i, \mathcal{O}_{(x,0), r_0})\}_i$ also forms a Cauchy sequence, the sets $\{L_{r_i} \cap [-1, 1]^2\}$ converge in the Hausdorff metric on compact subsets to the set $L_{r_{\infty}} \cap [-1, 1]^2$. And, since $0< \text{Height}(r_{\infty}, \mathcal{O}_{(x,0), r_0})\le 1/2$, we may estimate $\text{slope}( L_{r_{\infty}}) \in (0, 2/3]$.  Therefore, for any $0< \delta$ there exists an $i(\delta) \in \mathbb{N}$ such that 
\begin{align*}
    \dist_{\mathcal{H}}(L_{r_\infty} \cap [-1, 1]^2, L_{r_j} \cap [-1, 1]^2) \le \delta
\end{align*}
for all $j \ge i(\delta)$. Therefore, if $r' \in [2r_{i(\delta)}/3, r_{\infty}]$ and 
\begin{align*}
    \left(3\text{Height}(r_\infty, \mathcal{O}_{(x,0), r_0}) r' - 2 \text{Height}(r_{\infty}, \mathcal{O}_{(x,0), r_0}) \right) - \text{Height}(r', \mathcal{O}_{(x,0), r_0}) \ge 4\delta> 0,
\end{align*}
then 
\begin{align*}
    \left(3\text{Height}(r_j, \mathcal{O}_{(x,0), r_0}) r' - 2 \text{Height}(r_j, \mathcal{O}_{(x,0), r_0}) \right) - \text{Height}(r', \mathcal{O}_{(x,0), r_0}) \ge \delta> 0,
\end{align*}
and by the minimality in the definition of $r_{i+1}$, for all $j \ge i(\delta)$ $r' \not \in [2r_{j}/3, r_{\infty})$.  Letting $j \rightarrow \infty$ we may assume that $r' \not \in (2r_{\infty}/3, r_{\infty})$. Repeating the argument for $\delta \rightarrow 0$ shows that the claim holds.  This proves the lemma.
\end{proof}

Using orthogonal projection, we obtain the following simple corollary.

\begin{corollary}\label{integrals and projection}
Let $n=2$, $0< \gamma$, and let $u, (x, 0), \mathcal{O}, C_2, r_0, \rho$ be as in the statement of Lemma \ref{n=2 geometry}. Then
\begin{align*}
    \int_{\partial \mathcal{O}_{(x, 0), \rho} \cap [-1,1]^2} |x_2|^\gamma d\sigma \ge \int_{\mathcal{O}_{(x, 0), \rho} \cap \partial[-1,1]^2} \frac{1}{6C_2} |x_2|^\gamma d\sigma.
\end{align*}
\end{corollary}

\begin{proof}
By reflection, without loss of generality, we may assume that $\mathcal{O}_{(x,0), \rho} \cap [-1, 1]^2 \subset [0, 1]^2$. Define $\partial^+\mathcal{O}_{(x,0), \rho}$ to be the set
\begin{align*}
    \partial^+\mathcal{O}_{(x,0), \rho} := & \{ (x_1, x_2) \in \partial\mathcal{O}_{(x,0), \rho}: x_2= \text{Height}(|x_1|, \mathcal{O}_{(x,0), \rho}), \quad 2/3 \le x_1 \le 1\}\\
    & \cap \{(x_1, x_2) \in \partial\mathcal{O}_{(x,0), \rho}: x_2 \ge 3 \text{Height}(1, \mathcal{O}_{(x,0), \rho})x_1 - 2\text{Height}(1, \mathcal{O}_{(x,0), \rho}) \}.
\end{align*}
Let $\pi_1$ be orthogonal projection onto $\{x_2 = 0 \}$. Let $\pi_{2}$ be orthogonal projection onto the line $\{x_1 = 1\}$. If $f$ is the linear function such that $L_r = \text{graph}_{\mathbb{R}}(f)$, define $\pi_{L_r}: \mathbb{R}^2 \rightarrow L_r$ to be the function $\pi_{L_r}(x_1, x_2) = (x_1, f(x_1))$. Note that for $C_2 \le 1/2$, $|\nabla f| \le 3/2$. 

We observe that for any set $U\subset \mathbb{R}^2$, 
\begin{align*}
    \mathscr{H}^1(U) & \ge \mathscr{H}^1(\pi_1(U))\\
    & \ge \frac{1}{\sqrt{1+(3/2)^2}}\mathscr{H}^1(\pi_{L_r}(U)) > \frac{1}{2} \mathscr{H}^1(\pi_{L_r}(U)).
\end{align*}

Then (\ref{L_r}), Definition \ref{weak solution def}(3), and the choice of $C_2 \le \frac{1}{2}$ implies that 
\begin{align*}
    \int_{\partial \mathcal{O}_{(x, 0), \rho} \cap [-1,1]^2} |x_2|^\gamma d\sigma & \ge \int_{\partial^+\mathcal{O}_{(x, 0), \rho} \cap [-1,1]^2} |x_2|^\gamma d\sigma \\
    & \ge \frac{1}{2}\int_{\pi_{L_r} \left(\partial^+\mathcal{O}_{(x, 0), \rho}\right) \cap ([\frac{2}{3}, 1] \times [0, 1])} |x_2|^\gamma d\sigma\\
    & \ge \frac{1}{2}\int_{\pi_2(L_r \cap ( [\frac{2}{3}, 1] \times [0, 1]))} \sqrt{(3C_2)^{-2}+1}|x_2|^\gamma d\sigma\\
    & \ge \frac{1}{2}\int_{\mathcal{O}_{(x, 0), r_2} \cap (\{1\} \times [0, 1])} \frac{1}{3C_2}|x_2|^\gamma d\sigma.
\end{align*}
\end{proof}

\subsection{Proof of Theorem \ref{macroscopic information}}

Let $u$ satisfy (\ref{checkable}) with constant $C< \infty$ in $B_r((x, 0)) \subset \Omega$. We argue by contradiction.  Let $0< r_0$ be such that $0< \text{Height}(1, \mathcal{O}_{(x, 0), r_0}) \le  C_2 \le 1/2$.  Let $0< \rho \le r_0$ as in Lemma \ref{n=2 geometry}. We consider 
$$
V := \mathcal{O}_{(x,0), \rho} \cap \partial [-1, 1]^2.
$$
And, note that 
\begin{align}
\int_{V}|x_2|^{\gamma}d\sigma \le 
\int_{\partial D \cap [-1, 1]^2}|x_2|^\gamma d\sigma
\end{align}
for any set $D$ which is relatively open in $[-1, 1]^2$ and satisfies $V \subset (\partial [-1, 1]^2  \cap D)$.  

Next, use Definition \ref{weak solution def}(3), the Divergence Theorem, and (\ref{checkable}) to calculate 
\begin{align*}
    \int_{\partial \mathcal{O}_{(x,0), \rho} \cap [-1, 1]^2} |x_2|^\gamma d\sigma & \le \Delta u(\overline{\mathcal{O}_{(x,0), \rho}} \cap [-1, 1]^2)\\
    & = \int_{V} \nabla u \cdot \vec{\eta} d\sigma\\
    & \le \int_{V} C |x_2|^\gamma d\sigma.
\end{align*}

But, by Corollary \ref{integrals and projection}
\begin{align*}
    \int_{\partial \mathcal{O}_{(x,0), \rho} \cap [-1, 1]^2} |x_2|^\gamma d\sigma & \ge \int_{V_r} \frac{1}{6C_2} |x_2|^\gamma d\sigma.  
\end{align*}
Therefore, taking $6C_2 \le \frac{1}{C}$ we have a contradiction. This proves Theorem \ref{macroscopic information}.

\bibliography{references}
\bibliographystyle{amsalpha}

\end{document}